# Remarks on normal bases

Marcin Mazur

*Department of Mathematics, the University of Chicago,*

*5734 S. University Avenue, Chicago, IL 60637*

e-mail: mazur@math.uchicago.edu

**Abstract**

We prove that any Galois extension of commutative rings with normal basis and abelian Galois group of odd order has a self dual normal basis. Also we show that if $S/R$ is an unramified extension of number rings with Galois group of odd order and $R$ is totally real then the normal basis does not exist for $S/R$.

## 1  Introduction.

Let $R \subset S$ be an extension of commutative rings. Suppose that $G$ is a finite group of automorphisms of $S$ and $R = S^G$. By $S^{(G)}$ we denote the ring $Map(G, S)$ of all functions from $G$ to $S$. The group $G$ acts on $S^{(G)}$ by $f^g(h) = f(gh)$ and $S$ (= constant functions) is the ring of invariants. There is an obvious ring homomorphism $\phi : S \otimes_R S \longrightarrow S^{(G)}$ given by $\phi(s_1 \otimes s_2)(g) = s_1 s_2^g$. This map is $G$−equivariant, where $G$ acts on $S \otimes_R S$ via the second component.

Recall that $S$ is called Galois over $R$ if $\phi$ is surjective (and then $\phi$ is in fact an isomorphism). Clearly $S^{(G)}/S$ is Galois. If $R$, $S$ are Dedekind domains then the extension $S/R$ is Galois iff the corresponding extension of fields of fractions is Galois with group $G$, $S^G = R$ and $S/R$ is unramified (see [3] for more about Galois extensions of rings).

It is well known that if $S/R$ is Galois then $S$ is a projective, faithfully flat $R$−module of constant rank $|G|$. Moreover the trace map $tr : S \longrightarrow R$, $tr(s) = \sum s^g$ coincides with the module-theoretic trace and is surjective (for the proof base change to $S$ where the situation is clear and then use f.f. descent).

The classical problem is to describe the structure of $S$ as an $RG-$module. First step is the following well known lemma:

**Lemma 1** *Let $S/R$ be an extension of commutative rings such that $R = S^G$, $S$ is $R-$projective and the trace map is surjective. Then $S$ is a projective $RG-$module.*

*Proof:* There exists an $RG-$modules epimorphism $p : F \longrightarrow S$ with $F$ a free $RG-$module. Since $S$ is $R-$projective there exists an $R-$module splitting $f$ of $p$. Let $c \in S$ be such that $tr(c) = 1$. Define a new map $h : S \longrightarrow F$ by $h(s) = \sum g^{-1} f(s^g c)$. Clearly $h$ is an $RG-$module map and $ph(s) = \sum g^{-1}(s^g c) = str(c) = s$. $\square$

The natural question to ask is under what circumstances $S$ is a free $RG-$module. For rings of integers in finite extensions of the rationals (we call such rings *number rings*) this is an old and still unsolved problem.

If $S$ is a free $RG-$module then there is an element $s \in S$ such that the orbit of $s$ under $G$ is an $R-$basis. We call any such basis a *normal basis* of $S$. If moreover this basis is self-dual with respect to the trace form then we call it a *self-dual normal basis*. For example, the extension $S^{(G)}/S$ has always a self dual normal basis generated by $\delta_-^1$ where

$$\delta_-^1(g) = \delta_g^1 = \begin{cases} 1 & if \ g = 1 \\ 0 & if \ g \neq 1 \end{cases}$$

## 2 Self dual normal basis.

The following useful theorem should be well known to the experts but we do not know of any good reference:

**Theorem 1** *Let $S$ be a commutative ring, $G$ a finite group of automorphisms of $S$ and $R = S^G$. For $u \in S$ the following are equivalent:*
*1) $\sum u^g g^{-1}$ is a unit in $SG$;*
*2) $S/R$ is Galois with normal basis generated by $u$.*

*Proof:* Let $(\sum u^g g^{-1})(\sum w_g g) = 1$. In other words, for any $h \in G$ we have $\sum u^g w_{gh} = \delta_h^1$, i.e. $\sum u^{gh^{-1}} w_g = \delta_h^1$, so also $\sum u^g w_g^h = \delta_h^1$. Since $\phi(\sum u^g \otimes w_g)(h) =$ $\sum u^g w_g^h = \delta_h^1$, the natural map $\phi: S \otimes_R S \longrightarrow S^{(G)}$ is surjective and $S/R$ is Galois. In particular, $\phi$ is an isomorphism of $SG$−modules. Now $\phi(1 \otimes u) = (\sum u^g g^{-1})\delta_-^1$ is a free generator of the $SG$−module $S^{(G)}$, so $1 \otimes u$ is a free generator of $S \otimes_R S$. Since $RGu \subseteq S$ and after tensoring with $S$ we get equality, $u$ is a free generator of $S$ (note that $S/R$ is faithfully flat). This shows that 1) implies 2).

To get the converse observe that the $SG$−module isomorphism $\phi$ maps $1 \otimes u$ to a free generator $f: g \mapsto u^g$. But $\delta_-^1$ is also a free generator and $(\sum u^g g^{-1})\delta_-^1 = f$ so $\sum u^g g^{-1}$ is a unit. $\square$

**Corollary 1** *If $S/R$ is a Galois extension of commutative rings with group $G$ and normal basis then for any normal subgroup $H$ of $G$ the extension $S^H/R$ is Galois and has a normal basis.*

*Proof:* If $\sum u^g g^{-1}$ is a unit of $SG$ then under the natural surjection $SG \longrightarrow SG/H$ it maps to a unit $\sum v^h h^{-1}$ of $SG/H$, where $v = tr_{S/S^H} u$. But this is a unit in $S^H G/H$ so the result follows by Theorem 1. $\square$

Observe that if $(\sum u^g g^{-1})(\sum w_g g) = 1$ then also $(\sum u^{gh} g^{-1})(\sum w_g^h g) = 1$ so $(\sum u^{gh}(gh)^{-1})(\sum w_g^h gh) = 1$. Since the inverse is unique we get $w_g^h = w_{gh}$ for any $g, h \in G$. Putting $w = w_1$ we obtain $w_g = w^g$, so the inverse to $\sum u^g g^{-1}$ is $\sum w^g g$.

Recall that in a group ring $SG$ we have an involution $*$ which on $G$ acts as an inverse (i.e. $(\sum a_g g)^* = \sum a_g g^{-1}$). This shows that $\sum w^g g^{-1}$ is a unit and therefore $w$ generates a normal basis.

Suppose now that $u$ and $w$ generate a normal basis. Consider the unit $(\sum u^g g^{-1})(\sum w^g g^{-1}) = \sum c_g g^{-1}$ where $c_g = \sum u^h w^{gh^{-1}}$. We have $c_g^f = \sum u^{hf} w^{gh^{-1}f} = \sum u^h w^{g(hf^{-1})^{-1}f}$ for any $f \in G$. If $G$ is abelian then we get $c_g^f = \sum u^h w^{gh^{-1}f^2}$. In particular, for $c = c_1$ we get $c^f = \sum u^h w^{f^2 h^{-1}} = c_{f^2}$. If moreover the order of $G$ is odd then $\sum c_g g^{-1} = \sum c_{g^2} g^{-2} = \sum c^g g^{-2}$ so $\sum c^g g^{-2}$ is a unit. Since the map $g \mapsto g^2$ is an automorphism of $G$, also $\sum c^g g^{-1}$ is a unit. Thus $c$ generates a normal basis too. If $U = \sum u^g g^{-1}$, $W = \sum w^g g^{-1}$ and $C = \sum c^g g^{-1}$ then $C = \psi(UW)$, where $\psi$ is the automorphism of $SG$ which maps $g^2$ to $g$ for all $g \in G$. In particular, if $w$ is such that $UW^* = 1$ then $(UW)(UW)^* = 1$ and therefore $CC^* = 1$. In other words $c$ generates a self dual normal basis (with respect to the trace form). Therefore we proved the following

**Theorem 2** *If $S/R$ is a Galois extension of commutative rings with abelian Galois group of odd order and if it has a normal basis then it has a self dual normal basis.*

For cyclic groups of odd order this result has been obtained by Kersten and Michaliček [4]. Note also that Bayer and Lenstra ([1], [2]) proved that for fields and any group of odd order a self dual normal basis exists. The author does not know whether this remains true for Galois extensions of rings.

## 3 Number rings.

Suppose that $S/R$ is a Galois extension of number rings with an abelian Galois group $G$ of odd order and having a normal basis. By Theorem 2 it has a self dual normal basis generated by $a \in S$. Thus $X = \sum a^g g^{-1}$ is a unit in $SG$ and $XX^* = 1$. We assume that $R$ is totally real (and so is $S$). Consider any ring homomorphism $\psi : SG \longrightarrow \mathbb{C}$. Since $S$ is totally real, we have $\psi(u^*) = \overline{\psi(u)}$ for any $u \in SG$. In particular $\psi(X^2) = \psi(X/X^*) = \psi(X)/\overline{\psi(X)}$ is of absolute value 1. Note that $\psi(X^2)$ is an algebraic integer and we just proved that all of its conjugates have absolute value 1 ($\psi$ is arbitrary) so by a well known theorem of Kronecker $\psi(X^2)$ is a root of 1. Since this holds for all $\psi$, we conclude that $X^2$ is a torsion unit of $SG$ (note that $\mathbb{C}G = \mathbb{C}^{|G|}$) and so is $X$. Since no prime divisor of the order of $G$ is invertible in $S$ and $\pm 1$ are the only torsion units of $S$, the order of $X$ is a divisor of $2|G|$ (in general, if $G$ is a finite group and $S$ a commutative ring such that no prime divisor of the order of $G$ is invertible in $S$ then any unit of finite order in $RG$ is of the form $az$ where $a$ is a torsion unit of $S$ and the order of $z$ divides $|G|$; for a proof see for example [6], Lemma 5). Now note that the trace of a regular representation of SG is given by $T(\sum u_g g) = |G|u_1$. But the trace of an element of finite order dividing $2|G|$ is a sum of $2|G|$−th roots of 1 so in particular $|G|a = T(X) \in \mathbb{Q}(\xi_{|G|})$. Let $K$, $L$ be the fields of fractions of $R$, $S$ respectively. Thus we showed that $L = K(a) \subseteq K(\xi_{|G|})$. But this is clearly false since $[L : K] = |G| > [K(\xi_{|G|}) : K]$. The contradiction shows that $S/R$ can not have a normal basis. Therefore we proved the following surprising theorem:

**Theorem 3** *Let $R$ be the ring of integers in a totally real number field. If $S/R$ is a Galois extension of number rings with Galois group $G$ of odd order then it does not*

*have a normal basis.*

*Proof:* For abelian $G$ the result was shown above. In general $G$ is solvable so it has a non trivial abelian quotient $H$. If $S/R$ had a normal basis then so would $S^H/R$ by Corollary 1. But this is false since $S^H/R$ is Galois with abelian Galois group of odd order. $\square$

Note that Taylor ([7]) proved that in the situation of the above theorem $S$ is always a free $\mathbb{Z}G-$module.

The above method allows us to prove also the following

**Proposition 1** *Let $R = \mathbb{Z}[\xi_{p^k} + \xi_{p^k}^{-1}]$ where $p$ is an odd prime. Let $S$ be the ring of integers in a cyclic extension of $\mathbb{Q}(\xi_{p^k} + \xi_{p^k}^{-1})$ of degree $p^n$. If $S[1/p]/R[1/p]$ is Galois and has a normal basis then it coincides with the cyclotomic $p^n-$extension.*

*Proof:* We keep the above notation. First we show that $\psi(X^2)$ is an algebraic integer. Clearly $\psi(R) = R$, $\psi(SG) \subseteq \psi(S)[\xi_{p^n}]$ and $\psi(S)$ is the ring of integers in a cyclic extension $L$ of $\mathbb{Q}(\xi_{p^k} + \xi_{p^k}^{-1})$ of degree $p^n$, so we will write $S$ for $\psi(S)$. Note that all primes of $S[\xi_{p^n}]$ over $p$ are stable under complex conjugation. To show this let $m = \max\{n, k\}$. There is only one prime $\pi$ over $p$ in $\mathbb{Q}(\xi_{p^m} + \xi_{p^m}^{-1})$ and it ramifies in $\mathbb{Q}(\xi_{p^m})$. Thus all primes of $L(\xi_{p^m} + \xi_{p^m}^{-1})$ over $\pi$ ramify in $L(\xi_{p^m})$ and therefore are stable under complex conjugation (note that $L(\xi_{p^m})/\mathbb{Q}(\xi_{p^m})$ is of $p-$power degree). Now $\psi(X)$ is a $p-$unit in $L(\xi_{p^m})$ and therefore $\psi(X)/\overline{\psi(X)} = \psi(X^2)$ is an algebraic integer. As before we conclude that $\psi(X)$ is a root of 1 in $L(\xi_{p^m})$. Therefore $\psi(X)^{2p^m} = 1$ and consequently $X^{2p^m} = 1$. The argument with traces shows now that $L = K(a) \subseteq \mathbb{Q}(\xi_{p^m})$ and this implies that $L$ is the cyclotomic $p^n-$extension of $\mathbb{Q}(\xi_p + \xi_p^{-1})$. $\square$

**Remark.** The cyclotomic $p^n-$extension has a normal basis, as shown in [3].

As a direct consequence we get the following result of Kersten and Michaliček ([5]):

**Corollary 2** *If $k = n = 1$ and $S/R$ is Galois then $S[1/p]/R[1/p]$ does not have a normal basis.*

**Remark.** Corollary 2 suggests the following attack on Vandiver's Conjecture: show that any extension as above has to have a normal basis and derive that there is no such extensions. Of course at present nobody knows how to do that.

# References

[1] E. Bayer-Fluckiger and H. W. Lenstra, Jr., *Forms in odd degree extensions and self-dual normal bases*, Amer. J. Math. 112 (1990), 359–373.

[2] E. Bayer-Fluckiger, *Self-dual normal bases*, Indag. Math. 51 (1989), 379–383

[3] C. Greither, Cyclic Galois extensions of Commutative Rings, LNM 1534, Springer–Verlag, Berlin, 1992.

[4] I. Kersten and J. Michaliček, *Kubische Galoiserweiterungen mit Normalbasis*, Comm. in Algebra 9 (1981),1863–1871.

[5] I. Kersten and J. Michaliček, *A remark about Vandiver's Conjecture*, C.R. Math. Rep. Acad. Sci. Canada, vol. VII (1985), 33–37.

[6] M. Mazur, On the izomorphism problem for infinite group rings, Expo. Math 13 (1995), 433-445.

[7] M. J. Taylor, *On Frochlich conjecture for rings of integers of time extensions*, Invent. Math. 63 (1981), 41–79.